\documentclass{article}
\usepackage{amsthm}
\usepackage{amsmath,a4wide}
\usepackage{amssymb}
\usepackage{color}

\theoremstyle{theorem}
\newtheorem{theorem}{Theorem}

\theoremstyle{definition}

\newtheorem*{example}{Example}

\newtheorem{prop}{Proposition}

\DeclareMathOperator{\p}{per}
\DeclareMathOperator{\pp}{pp}

\begin{document}
\title{Perfect Overpartitions and Factorization  of Integers}  
\author{Augustine O. Munagi\\
{\small School of Mathematics, University of the Witwatersrand,}\\
{\small P O Wits 2050, Johannesburg, South Africa.}\\ 
{\small Augustine.Munagi@wits.ac.za} }
\date{}
\maketitle
 
\begin{abstract}
In his classic text, \emph{Combinatory Analysis}, MacMahon defined a perfect partition of a positive integer $n$ as a partition whose parts contain exactly one partition of every positive integer not exceeding $n$. In this paper we apply the same definition to overpartitions which are integer partitions with the additional property that the final occurrence of each part may be overlined. It turns out that perfect overpartitions are enumerated by ordered factorization functions in which the occurrence of 2 as a factor determines the presence of an overlined part. 
\medskip

\smallskip

\end{abstract}

\section{Introduction}\label{introd}
A partition of a positive integer $n$ is any weakly increasing sequence of positive integers that sum to $n$. The summands are called \emph{parts}, and $n$ is the \emph{weight}, of the partition. A partition $\lambda$ of $n$ (often expressed as $\lambda\vdash n$) into $k$ parts will be denoted by 

\begin{equation*}
\lambda = (\lambda_1,\lambda_2,\ldots,\lambda_k),\, 1\leq\lambda_1\leq\cdots\leq\lambda_k>0
\end{equation*}
or 
\begin{equation*}
\lambda = (\lambda_1^{m_1},\lambda_2^{m_2},\ldots,\lambda_r^{m_r}),\, 0<\lambda_1<\cdots < \lambda_r,\, 1\leq r\leq k,
\end{equation*}
where $m_i$ denotes the multiplicity of $\lambda_i$ for all $i$.

A \emph{perfect partition} of $n$ is a partition in which the parts contain exactly one partition of every positive integer less than or equal to $n$. For example, $(1^3,4)\vdash 7$ is a perfect partition since it contains the partitions $(1)$, $(1^2)$, $(1^3)$, $(4),(1,4),(1^2,4),(1^3,4)$ with weights $1,2,\ldots,7$, respectively. 

There is a known bijection between the set of perfect partitions of $n$ and the set of {ordered factorizations} of $N=n+1$, namely $a_1a_2\cdots a_r,\, a_i>1$ for all $i$ given by (see \cite{An,Mac2,R}) 
\begin{equation} \label{eqbij}
a_1 a_2\cdots a_r \longmapsto \left(1^{a_1-1},a_1^{a_2-1},(a_1 a_2)^{a_3-1},\ldots, (a_1 a_2\cdots a_{r-1})^{a_r-1}\right).
\end{equation}

Let $f(n,k)$ be the number of ordered factorizations of $n$ into $k$ factors, and let the prime-power factorization of $n$ be $n=p_1^{\alpha_1}p_2^{\alpha_2}\cdots p_r^{\alpha_r}$.
The formula for $f(n,k)$ is given by (see \cite{Mac2} or \cite[p.\ 59]{An})
\begin{equation} \label{eqfnk}
f(n,k)=\sum_{i=0}^{k-1}(-1)^{i}\binom{k}{i} \prod_{j=1}^{r} \binom{\alpha_j+k-i-1}{\alpha_j},\ n>1,
\end{equation}
where we adopt the convention $f(1,0)=1$. Then we define 
\begin{equation} \label{eqfn}
f(n):= f(n,1)+f(n,2)+\cdots + f(n,\Omega(n))
\end{equation}
with $f(n)=0,\,n<1$ and $f(1)=1$, where $\Omega(n)$ is the number of prime factors of $n$ counted with multiplicity.

In many cases it is more efficient to find $f(n)$ using the following classical recurrence relation (see, for example, \cite{KM}).  
\begin{equation} \label{off1}
f(1)=1,\quad f(n)=\sum_{\scriptstyle d\mid n\atop \scriptstyle d<n}f(d).
\end{equation}

From the bijection \eqref{eqbij} the formula for the number $\p(n)$ of perfect partitions of $n$ is given by
\begin{equation} \label{eqppform}
\p(n) = f(n+1)
\end{equation}
while $f(n,k)$ counts perfect partitions of $n$ with $k$ blocks or runs of equal parts.

\begin{example}
Table \ref{tabperf0} shows the ordered factorizations of 6 which correspond to the perfect partitions of 5.

\begin{table}[h]
\begin{center}
\begin{tabular}{l|c|c|c} \hline
Ordered Factorization of 6 &6 &$2\cdot 3$ &$3\cdot 2$\\ \hline
Perfect Partition of 5 &$(1^5)$ &$(1,2^2)$ &$(1^2,3)$\\ \hline
\end{tabular}
\caption{Factorizations of 6 and perfect partitions of 5}\label{tabperf0}
\end{center}
\end{table}

\end{example}

The definition of a perfect partition first appeared in the works of P. A. MacMahon \cite{Mac1, Mac2}. Subsequently other mathematicians studied and found several properties and generalizations of perfect partitions (see for example, \cite{AS1, AS2, Lee, Mu1, Par1}). 
Park \cite{Par1} generalized perfect partitions to `complete partitions' by removing the uniqueness condition from subpartitions, and considered partitions of $n$ that contain at least one partition of every positive integer $<n$.
Agarwal and Sachdeva \cite{AS1} studied properties of $n$-color perfect partitions that contain one $n$-color partition of every positive integer $<n$. More recently Munagi \cite{Mu1} explored perfect compositions defined analogously.

The purpose of this paper is to extend the unique representation property of perfect partitions to overpartitions. 
An overpartition of $n$ is any partition of $n$ with the additional property that the last occurrence of each part may be distinguished by being overlined.   

For example, the overpartitions of 3 are 
$$(3), (\overline{3}), (1,2), (1,\overline{2}), (\overline{1},2), (\overline{1}, \overline{2}),
(1,1,1), (1,1,\overline{1}).$$

In Section \ref{overpartn} we define perfect overpartitions. Then in Section \ref{formula} we explore their enumeration formulas. Finally, in Section \ref{popform} we derive exact formulas for numbers of perfect overpartitions of some positive integers with small 2-adic order.

\section{Perfect Overpartitions}\label{overpartn}
A \emph{perfect overpartition} of $n$ is any overpartition of $n$ that contains exactly one overpartition of every smaller positive integer.
A perfect partition is a perfect overpartition (with no overlined part).
An overlined part $\lambda_i$ in a perfect overpartition $\lambda$ is the only part of size $\lambda_i$ in $\lambda$.

A perfect partition $\lambda$ may generate a perfect overpartition by overlining the last occurrence of some parts in $\lambda$. The resulting overpartition has an overlined part only if $\lambda$ contains a distinct part-size. If it contains $r\geq 1$ distinct part-sizes, then it gives a perfect overpartition with at most $r$ overlined parts. If $r=0$, then $\lambda$ may still generate several overpartitions. However, none of the latter can be perfect since each would contain two overpartitions of some integer $a$, namely $(a)$ and  $(\overline{a})$. 

The bijection \eqref{eqbij} shows that any factor $a_i=2$ gives rise to a unique part-size in the corresponding perfect partition. Consequently, a perfect overpartition of $n$ containing an overlined part corresponds to an ordered factorization of $n+1$ containing 2.
So the weight of a perfect overpartition containing an overlined part is odd.

Let $\overline{\pp}(n)$ be the number of perfect overpartitions of $n$. 
If $n$ is even, then $\overline{\pp}(n)=\p(n)$ since no factorization of $n+1$ can contain 2. 
If $n$ is odd, then $\overline{\pp}(n)>\p(n)$ and $\overline{\pp}(n)-\p(n)$ enumerates objects with at least one overlined part.

Thus, for example, the factorization $3\cdot 2\cdot 2$ maps to the perfect partition $(1^2,3,6)\vdash 11$ which generates the four perfect overpartitions $(1^2,3,6), (1^2,\overline{3},6), (1^2,3,\overline{6}), (1^2,\overline{3},\overline{6})$. Each object contains a unique overpartition of $m=1,\ldots,11$. For instance, $(1^2,\overline{3},6)$ contains the overpartitions 
$$(1),(1^2),(\overline{3}),(1,\overline{3}),(1^2,\overline{3}),(6),(1,6),(1^2,6),(\overline{3},6),(1,\overline{3},6),(1^2,\overline{3},6).$$
Note that the further derived overpartition $(1,\overline{1},3,6)$ is not perfect because it contains two overpartitions of 1.

\section{Enumeration Formulas}\label{formula}
Let $\overline{\pp}(n,r)$ denote the number of perfect overpartitions of $n$ with $r$ overlined parts, $r\geq 0$, so that $\sum_r \overline{\pp}(n,r)=\overline{\pp}(n)$, and let $f_v(N)$ denote the number of ordered factorizations of $N$ containing $v$ copies of 2.

\begin{theorem}\label{thm1}
We have the following: 
\begin{equation}\label{eqn1}
\overline{\pp}(n,r)=\sum_{v=r}^s\binom{v}{r}f_{v}(n+1),\ r\geq 0,
\end{equation}
where $n+1=2^s m,\, 2\nmid m, s\geq 0$. Hence
\begin{equation}\label{eqn2}
\overline{\pp}(n)=\sum_{r=0}^s\sum_{v=r}^s\binom{v}{r}f_v(n+1).
\end{equation}
\end{theorem}
In particular the first summand gives $\, \sum_{v=0}^s f_v(n+1) = f(n+1) = \p(n)$. 
\begin{proof}[Proof of Theorem \ref{thm1}] Let $R$ be an ordered factorization of $n+1$ containing $v$ copies of 2, $0\leq v\leq s$, and let $\beta$ be the perfect partition corresponding to $R$ by the bijection \eqref{eqbij}. Then $\beta$ contains exactly $v$ distinct part-sizes. We may overline any $r$ of the latter, and there are $\binom{v}{r}$ choices of doing so. Since there are $f_v(n+1)$ possible factorizations $R$, the result follows.  
\end{proof}
\vskip 3pt

\begin{example} We verify that $\overline{\pp}(11)=19$ by means of the bijection \eqref{eqbij} and overlining of distinct parts. The ordered factorizations of $12=2^2 3$ (with $s=2$) are as follows:
\begin{equation}\label{eqfact12}
12,\, 2\cdot 6,\, 6\cdot2,\, 3\cdot4,\, 4\cdot3,\, 2\cdot2\cdot3,\, 2\cdot3\cdot2,\, 3\cdot2\cdot2.
\end{equation}
Thus $f_0(12)=3,f_1(12)=2,f_2(12)=3$, and (with $f_v:=f_v(12)$) we obtain the following: 
\medskip

All eight factorizations give $\, \overline{\pp}(11,0)=f_0+f_1+f_2=8$:
\[(1^{11}),(1,2^5),(1^2,3^3),(1^3,4^2),(1,2,4^2),(1^5,6),(1,2^2,6),(1^2,3,6).\]
 
The five factorizations containing 2 give $\, \overline{\pp}(11,1)=f_1+2f_2=8$:    
\[(\overline{1},2^5),(\overline{1},2,4^2),(1,\overline{2},4^2),(1^5,\overline{6}),(\overline{1},2^2,6),(1,2^2,\overline{6}),(1^2,\overline{3},6),(1^2,3,\overline{6}).\]
\vskip 5pt

The three factorizations with two 2's give $\, \overline{\pp}(11,2)=f_2=3$: 
\[(\overline{1},\overline{2},4^2),(\overline{1},2^2,\overline{6}),(1^2,\overline{3},\overline{6}).\]
Therefore, $\, \overline{\pp}(11) = \overline{\pp}(11,0) + \overline{\pp}(11,1) + \overline{\pp}(11,2) = 19$.

\end{example}

It is clear that effective computation of $\overline{\pp}(n)$ for any $n>0$ relies on a formula for $f_v(n+1)$. Note that $n+1=2^s m,\, s\geq 0$, where $m$ is odd. So $s$ is the 2-adic order of $N=n+1$ which is usually denoted by $\nu_2(N)$. It's exact value may be obtained using  Legendre's theorem for the highest power of a prime dividing a factorial:

\begin{equation}\label{eqn2adic}
s=\nu_2(N)=\sum_{i\geq 1}\left(\left\lfloor\frac{N}{2^i}\right\rfloor - \left\lfloor\frac{N-1}{2^i}\right\rfloor\right).
\end{equation}

The computations in this paper were performed by using the computer algebra system \emph{Maple}.
\vskip 4pt

Proofs of some theorems below will use the following results on compositions or ordered partitions (see, for example, \cite{An,HM}).
\medskip

The number of compositions of $n$ into $k$ parts is given by 
\begin{equation}\label{eqcompk1}
c(n,k)=\binom{n-1}{k-1}.
\end{equation}

The number of compositions of $n$ without 1's into $k$ parts is given by 
\begin{equation}\label{eqcompk2}
c_2(n,k)=\binom{n-k-1}{k-1}.
\end{equation}

We will also apply a special case of the Vandermonde's convolution identity (see \cite{R}):
\begin{equation}\label{eqvand}  
\sum_{i\geq 0}\binom{X}{i}\binom{Y}{Y-i} = \binom{X+Y}{Y}
\end{equation}
for any nonnegative integers $X,Y$.

We will sometimes denote by $\{r(n)\}$ the set which is enumerated by the function $r(n)$. 

\begin{prop}\label{prop1}
Let $N=2^s m$ with $m$ odd, $m>1$. Then the number $g_s(N)$ of ordered factorizations of $N$ containing $s$ copies of 2 is given by   
\begin{equation}\label{eqn4b}
g_s(N)=\sum_{j=1}^{\Omega(m)}\binom{s+j}{s}f\left(m,j\right).
\end{equation}
\end{prop}
\begin{proof} We construct an element $R\in \{g_s(N)\}$. If $s=0$, then $R$ consists of odd factors: $g_0(N)=\sum_{j\geq 1} f\left(m,j\right) = f(m)$.

So we assume that $s>0$ with $m>1$ and obtain $R$ with $s+j$ ($j>0$) factors as follows. Designate $s+j$ cells for the sequence of factors and select any $j$ cells (in $\binom{s+j}{j}$ ways) to hold an ordered factorization of $m$ into $j$ factors (with $f(m,j)$ possible objects). Lastly, fill the remaining cells with $s$ copies of 2 in one way. Thus altogether we obtain $\binom{s+j}{j}f(m,j)$ objects $R$ for each $j$. Hence the result follows.
\end{proof}

By slightly modifying the proof of Proposition \ref{prop1} we obtain $f_{s-1}$ and $f_{s-2}$.

\begin{prop}\label{prop2} Let $N=2^s m$ with $m$ odd, $m>1$. Then we have the following.
\begin{align}
\text{(i)}\ f_{s-1}(N)& = \sum_{j\geq1}j\binom{j+s-1}{s-1}f(m,j)\label{eqsm1},\ s\geq 1,\\
\text{(ii)}\ f_{s-2}(N)& = \sum_{j\geq1}\left[(s-1)\binom{j+s-1}{s-1}+\binom{j+1}{2}\binom{j+s-2}{s-2}\right ]f(m,j),\ s\geq 2.\label{eqsm2}
\end{align}
\end{prop} 
\begin{proof} (i) If $s-1$ copies of 2 are factors, we multiply the remaining 2 into any of $j$ odd factors.\\
(ii) When 4 and $s-2$ copies of 2 are factors, we first put 4 into any of $s-1$ possible cells before the 2's are inserted uniquely. When 4 is not a factor, we multiply the remaining $2^2$ into any of $j$ odd factors taking $i$ at a time via $c(2,i),\, i\geq 1$. Thus altogether $f_{s-2}$ is given by 
$$\sum_{j\geq1}\binom{j+1+(s-2)}{j}f(m,j)(s-1) + \sum_{j\geq1}\binom{j+s-2}{j}f(m,j)\sum_{i\geq1}\binom{j}{i}c(2,i).$$
Lastly note that $\sum_{i\geq1}\binom{j}{i}c(2,i) = \binom{j+1}{2}$ from \eqref{eqcompk1} and \eqref{eqvand}.
\end{proof}
\medskip

A general formula can be stated for $f_v(2^s), 0\leq v<s$ with the observation that $f_s(2^s)=1$. A factorization is obtained by first designating $v+i$ cells to hold $v$ copies of 2 and $i$ higher powers of 2, $i>0$. Then there are $\binom{v+i}{i}$ ways to fix $i$ cells for higher powers. The  remaining exponent $s-v$ is then distributed among the $i$ cells using \eqref{eqcompk2}, whence each composition $(s_1,\ldots,s_i),\, s_j>1$ corresponds to $(2^{s_1},\ldots,2^{s_i})$. This gives as many factorizations as
$$\binom{v+i}{i}c_2(s-v,i),\ i\geq 1.$$
Thus we have proved the following assertion.
\begin{prop}\label{prop3} Let $s,v$ be integers with $0\leq v\leq s$. The number of ordered factorizations of $2^s$ which contain $v$ copies of 2 is given by
\begin{equation}\label{eqpow2}
f_v(2^s) = \delta_{vs} + \sum_{i=1}^{\lfloor\frac{s-v}{2}\rfloor}\binom{v+i}{v}\binom{s-v-i-1}{i-1},
\end{equation}
where $\delta_{ij}$ is the Kronecker delta ($i=j\implies \delta_{ij}=1$ and $i\neq j\implies \delta_{ij}=0$).
\end{prop} 
\medskip

At the other extreme are odd integers $N>1$. Since a factorization of $N$ cannot have a 2, Equation \eqref{eqfn} or \eqref{off1} immediately implies the following:
\begin{equation}\label{eqoddf2}
N\equiv 1\ \text{(mod 2)}\ \implies f_v(N) = f(N)\delta_{v0}.
\end{equation}

Now let $N>0$ be an even non-power of 2. The next three theorems give enumeration formulas for cardinalities of the following mutually disjoint sets whose union is $\{f_v(N)\}$.
\begin{enumerate}
\item $\{f_v(N)_1\}=$ set of ordered factorizations of $N$ which contain $v$ copies of 2 and avoid higher powers of 2.
\item $\{f_v(N)_2\} =$ set of ordered factorizations of $N$ which contain $v$ copies of 2 with some higher powers of 2 and avoid even non-powers of 2.
\item $\{f_v(N)_3\} =$ set of ordered factorizations of $N$ which contain $v$ copies of 2, some higher powers of 2 and some even non-powers of 2.
\end{enumerate}

\bigskip

\begin{theorem}\label{thm2}
Let $N=2^s m,\, s\geq 1$ such that $m$ is odd, $m>1$, and let $0\leq v<s$. The number of ordered factorizations of $N$ which contain $v$ copies of 2 and avoid higher powers of 2 is given by
\begin{equation}\label{eqn4b1}
f_{v}(N)_1 = \sum_{j=1}^{\Omega(m)}\binom{v+j}{v}\binom{j+s-v-1}{s-v}f\left(m,j\right).
\end{equation}
\end{theorem}
\begin{proof} We construct a member $R\in \{f_{v}(N)_1\}$. First consider the enumeration function $g_v(N/2^{s-v})$. From the proof of \eqref{eqn4b} we obtain that the number of factorizations of $N/2^{s-v}$ into $j$ factors containing $v$ copies of 2 and a sequence of $j$ (odd) factors of $m=\frac{N}{2^s}$ is 
\begin{equation}\label{eqn4c}
\binom{v+j}{j}f\left(m,j\right).
\end{equation}
Now multiply some or all of the $j$ factors with the remaining powers of $2$ as follows. 
Select any $i$ from the $j$ factors, say $r_1,\ldots,r_i$, and obtain a composition of $s-v$ into $i$ parts, $(c_1,\ldots,c_i)$. Then multiply each $r_t$ with $2^{c_t}$ to obtain $R$.
Since there are $\binom{j}{i}$ possible combinations of the $j$ factors and $c(s-v,i)$ compositions, we obtain the factor
\begin{equation}\label{eqn4d}
\sum_{i\geq 1}\binom{j}{i}c(s-v,i) = \sum_{i\geq 1}\binom{j}{i}\binom{s-v-1}{i-1} = \binom{j+s-v-1}{s-v},
\end{equation}
where the last two equalities follow from \eqref{eqcompk1} and \eqref{eqvand} with the identity $\binom{n}{k}=\binom{n}{n-k}$. 
Thus for each $j\geq 1$ we obtain the following from \eqref{eqn4c} and \eqref{eqn4d} which gives the stated result:
\begin{equation}\label{eqn4e}
\binom{v+j}{j}f\left(m,j\right)\binom{j+s-v-1}{s-v}.
\end{equation}

\end{proof}
\vskip 4pt

\begin{theorem}\label{thm3}
Let $N=2^s m,\, s\geq 1$ such that $m$ is odd, $m>1$, and let $0\leq v\leq s-2$. The number of ordered factorizations of $N$ which contain $v$ copies of 2 with some higher powers of 2 and avoid even non-powers of 2 is given by 
\begin{equation}\label{eqn4b1}
f_{v}(N)_2 = \sum_{j=1}^{\Omega(m)}\sum_{i=1}^{\lfloor\frac{s-v}{2}\rfloor}\binom{v+i+j}{j}\binom{v+i}{v}\binom{s-v-i-1}{i-1}f(m,j).
\end{equation}
\end{theorem}
\begin{proof} Let $R\in \{f_{v}(N)_2\}$.
Designate $v+i+j$ to hold $v$ 2's, $i\geq 1$ higher powers of 2 and $j\geq 1$ odd factors. Then choose $j$ cells to hold odd factors and factorize $m$ into  corresponding $j$ factors, in as many ways as $\binom{v+i+j}{j}f(m,j), 1\leq j\leq \Omega(m)$. Then from the remaining $v+i$ cells choose $i$ to hold higher powers of 2 and obtain a composition of $s-v$ into corresponding $i$ parts $\geq 2$, in as many ways as $\binom{v+i}{i}c_2(s-v,i)$.\\
Thus for fixed $j$ and $i$ the total number of objects $R$ obtained is 
$$\binom{v+i+j}{j}f(m,j)\binom{v+i}{i}c_2(s-v,i).$$ 
Hence the result follows on applying \eqref{eqcompk2} and summing over $j$ and $i$.
\end{proof}
\medskip

When both higher powers of 2 and even non-powers of 2 are included in every factorization we obtain $f_{v}(N)_3$.
\begin{theorem}\label{thm4}
Let $N=2^s m,\, s\geq 1$ such that $m$ is odd, $m>1$, and let $0\leq v\leq s-3$. Then the number of ordered factorizations of $N$ which contain $v$ two's with some higher powers of 2 and some even non-powers of 2 is given by 
\begin{equation}\label{eqn4c1}
f_{v}(N)_3 = \sum_{j=1}^{\Omega(m)}\sum_{r=2}^{s-v-1}\sum_{i=1}^{\lfloor\frac{s-v}{2}\rfloor}\binom{v+i+j}{j}f(m,j)\binom{v+i}{i}\binom{r-i-1}{i-1}\binom{j+s-v-r-1}{s-v-r}.
\end{equation}
\end{theorem}
\begin{proof} We combine the proofs of Theorems \ref{thm2} and \ref{thm3} to ensure that each factorization simultaneously includes higher powers of 2 and even non-powers of 2, besides $v$ copies of 2. Note that $s-v$ is now split into two parts $r$ and $s-v-r$, where $r$ forms factors that are higher powers of 2 via $c_2(r,i)$ while $s-v-r$ determines even non-powers of 2 via $c(s-v-r,t)$, as described in the proofs of Theorems \ref{thm2} and \ref{thm3}:
$$\binom{v+i+j}{j}f(m,j)\binom{v+i}{i}c_2(r,i)\sum_{t\geq 1}\binom{j}{t}c(s-v-r,t),\ j\geq 1,\, r\geq 2,\, i\geq 1.$$
Then, on invoking Equations \eqref{eqcompk1}, \eqref{eqcompk2} and \eqref{eqvand} we obtain
$$\binom{v+i+j}{j}f(m,j)\binom{v+i}{i}\binom{r-i-1}{i-1}\binom{j+s-v-r-1}{s-v-r}.$$
Hence the result follows.
\end{proof}

We can now state the full formula in terms of the foregoing ones.

\begin{theorem}\label{thm5}
Let $N=2^s m,\, s\geq 1$ such that $m$ is odd, $m>1$. The number of ordered factorizations of $N$ containing $v$ copies of 2, $0\leq v\leq s-1$, is given by the formula
\begin{equation}\label{eqn4b20}
f_v(N) = f_v(N)_1 + f_v(N)_2 + f_v(N)_3,
\end{equation}
where $f_v(N)_1, f_v(N)_2$ and $f_v(N)_3$ are given in Theorems \ref{thm2}, \ref{thm3} and \ref{thm4} respectively.
\end{theorem}
\medskip

One may compute values of $\overline{\pp}(n)$ using Theorem \ref{thm1} in conjunction with Theorem \ref{thm5}. The formulas simplify when $m$ assumes certain special forms.

For instance, if $m=p^{\alpha}$, where $p$ is an odd prime, then $f(p^{\alpha},j)$ equals the number of compositions of $\alpha$ into $j$ parts: 
$f(p^{\alpha},j) = c(\alpha,j) = \binom{\alpha-1}{j-1}$. If $m=p_1p_2\cdots p_t$, where the $p_i$ are distinct odd primes, then $f(p_1 p_2\cdots p_t,j) = j!S(t,j)$, where $S(n,k)$ is the Stirling number of the second kind (see for example \cite{Mu0}).

\begin{example} We illustrate Theorem \ref{thm5} by computing the exact value of $\overline{\pp}(479)$. Since $480=2^5(15)$ we first use the theorem to obtain $f_v(480)$ for $0\leq v\leq 5$ (see Table \ref{tab479}). Then we invoke Equation \eqref{eqn1} and obtain $\overline{\pp}(479,r)$ for $0\leq r\leq 5$, leading to the total value $\overline{\pp}(479)=5898$.
 
\begin{table}[h]
\centering
\small
\begin{tabular}[t]{|c|c|c|c|c|}\hline
$v$&$f_v(480)_1$&$f_v(480)_2$&$f_v(480)_3$&$f_v(480)$\\ \hline
0&13&38&87&138\\ \hline
1&32&102&132&266\\ \hline
2&51&72&132&255\\ \hline
3&64&140&0&204\\ \hline
4&65&0&0&65\\ \hline
5&48&0&0&48\\ \hline
$f(480)$&-&-&-&976\\ \hline
\end{tabular}
$\ $
\begin{tabular}[t]{|c|c|}\hline
$r$&$\overline{\pp}(479,r)$\\ \hline
0&976\\ \hline
1&1888\\ \hline
2&1737\\ \hline
3&944\\ \hline
4&305\\ \hline
5&48\\ \hline
$\overline{\pp}(479)$&5898\\ \hline
\end{tabular}
\caption{\small Calculation of $\overline{\pp}(479,r),\, 0\leq r\leq 5$ and $\overline{\pp}(479)$}\label{tab479}
\end{table}

\end{example} 
\medskip

We conclude this section by stating recursive formulas for $f_v$ in analogy to \eqref{off1}. In view of Equation \eqref{eqoddf2} it suffices to state such formulas for only even integers. 

\begin{theorem} \label{thm21} Let $N>0$ be an even integer.\\ 
If $N$ is a power of 2, then
\begin{equation}\label{off1x} 
f_v(2^s) = \delta_{vs} + \sum_{i=1}^{\lfloor\frac{s-v}{2}\rfloor}\binom{v+i}{v}\binom{s-v-i-1}{i-1},\quad v\geq 0.
\end{equation} 
If $N$ is not a power of 2, then
\begin{equation} \label{off3}
f_v(N) = f_{v-1}(N/2) + \sum_{\scriptstyle d\mid N\atop \scriptstyle 2d<N}f_v(d),\ v>0,
\end{equation}
and 
\begin{equation} \label{off2}
f_0(N) = 1+\sum_{\scriptstyle d>1,\, d\mid n\atop \scriptstyle 2d<N}f_0(d).
\end{equation}
\end{theorem}
\begin{proof} Equation \eqref{off1x} is precisely Proposition \ref{prop3} and provides `boundary conditions' for the two recurrences.
 
In \eqref{off3} each $R\in\{f_v(N)\},\, v>0$ is obtained as follows. If $d|N$ and $1<d<N/2$, then $R$ is obtained by appending $N/d\neq 2$ to each member of $\{f_v(d)\}$. The case $d=N/2$ introduces an additional 2 which is then appended to a  member of $\{f_{v-1}(N/2)\}$ to give $R$.

The proof is the same for \eqref{off2}, when $v=0$, except that 1 is added to account for the single-term factorization $(N)$.
\end{proof}

Note that in \eqref{off3} it suffices to use divisors $d$ with $\nu_2(d)\geq v$; other divisors contribute 0 to the sum.

\section{Some Perfect Overpartition Formulas}\label{popform}
If $s=0$, then $n=N-1$ is even and we obtain $\, \overline{\pp}(n) = \text{per}(n) = f(n+1)$. 

If $s=1$, then we obtain, from Propositions \ref{prop2}(i) and \ref{prop1}, that 
\begin{align*}  
& f_0 = \sum_{j\geq 1}jf(m,j);\\ 
& f_1 =\sum_{j\geq 1}\binom{j+1}{1}f\left(m,j\right),
\end{align*}
where $m=\frac{n+1}{2^s}\equiv \frac{N}{2^s}$.
Thus when $n=2m-1, m>1$ we obtain
\begin{align}  
\overline{\pp}(n) & = (f_0+f_1)+f_1 = f(n+1)+f_1\notag\\ 
& = f(n+1)+\sum_{j\geq 1}(j+1)f\left(m,j\right).\label{popform1}
\end{align}
We also note the following identity for any integer $N$ with $\nu_2(N)=1$:
\begin{equation}\label{eqn3id1}
f_1-f_0 = f(N/2).
\end{equation}

If $s=2$, we obtain, from Propositions \ref{prop2} and \ref{prop1}, that 
\begin{align*}
& f_0 = \sum_{j\geq 1}\binom{j+2}{2}f(m,j).\\
& f_1 = \sum_{j\geq1}j(j+1)f(m,j) = 2\sum_{j\geq1}\binom{j+1}{2}f(m,j)\\
& f_2 = \sum_{j\geq1}\binom{j+2}{2}f(m,j).
\end{align*}
From the formulas for $f_0$ and $f_2$ we obtain the following identity whenever $\nu_2(N)=2$:
\begin{equation}\label{eqn3id3}
f_0 = f_2.
\end{equation}
Thus
\begin{equation}\label{eqn3id4}
f(N) = 2f_0 + f_1\quad \text{or}\quad f(N) = f_1 + 2f_2.
\end{equation}

Using \eqref{eqn3id3} we obtain a simple formula for $\overline{\pp}(n)$ (with $n=2^2m-1, m>1$) as follows:
\begin{align}
\overline{\pp}(n) & = (f_0+f_1+f_2)+(f_1+2f_2)+f_2\notag\\
& = 2f(n+1)+f_0\notag\\
& = 2f(n+1) + \sum_{j\geq1}\binom{j+2}{2}f(m,j).\label{popform2}
\end{align}

If $s=3$, we use Theorems \ref{thm2} to \ref{thm4} with Theorem \ref{thm5} to obtain the formula
\begin{equation*}
f_0 = \sum_{j=1}^{\Omega(m)}\left[(j+1)^2+\binom{j+2}{3}\right]f(m,j).
\end{equation*}
Then from Propositions \ref{prop2} and \ref{prop1} we get
\begin{align*}
& f_1 = \sum_{j\geq1}\left[2\binom{j+2}{2} + (j+1)\binom{j+1}{2}\right]f(m,j),\\
& f_{2} = \sum_{j=1}^{\Omega(m)}j\binom{j+2}{2}f\left(m,j\right),\\
& f_{3} = \sum_{j=1}^{\Omega(m)}\binom{j+3}{3}f\left(m,j\right).
\end{align*}
Finally we obtain, for $n=2^3m-1,m>1$, 
\begin{align}
\overline{\pp}(n) & = (f_0+f_1+f_2+f_3)+(f_1+2f_2+3f_3)+(f_2+3f_3)+f_3\notag\\
& = 4f(n+1)-3f_0-2f_1+4f_3\notag\\
& = 4f(n+1)+\sum_{j\geq1}\left[\binom{j+2}{3}-(j+1)^2(j+3)\right]f(m,j).\label{popform3}
\end{align}

Such formulas emerge rapidly when $n+1$ is a power of 2. Note that Equations \eqref{eqpow2} and \eqref{eqn1} imply the following.
\begin{align}
\overline{\pp}(2^s-1,r) & =\sum_{v=r}^s\binom{v}{r}f_{v}(2^s)\notag\\
%& =\sum_{v=r}^s\binom{v}{r}\left[\delta_{vs} + \sum_{i\geq 1}\binom{v+i}{v}\binom{s-v-i-1}{i-1}\right ]\\
& = \binom{s}{r} + \sum_{v=r}^s\binom{v}{r}\sum_{i\geq 1}\binom{v+i}{v}\binom{s-v-i-1}{i-1}.\label{eqpop2r}
\end{align}
Hence 
\begin{equation}\label{eqpop2}
\overline{\pp}(2^s-1) = 2^s + \sum_{r=0}^s\sum_{v=r}^s\binom{v}{r}\sum_{i=1}^{\lfloor\frac{s-v}{2}\rfloor}\binom{v+i}{v}\binom{s-v-i-1}{i-1}.
\end{equation}
\medskip

The sequence $\overline{\pp}(n), n\geq 1$ begins as follows (not yet in the database \cite{Sl}):
$$2, 1, 5, 1, 5, 1, 13, 2, 5, 1, 19, 1, 5, 3, 34, 1, 13, 1, 19, 3, 5, 1, 65, 2, 5,\ldots $$

Expectedly, the sequence for even weights, namely $\overline{\pp}(2n), n\geq 1$, is in \cite[A352063]{Sl} under the description, ``Number of ordered factorizations of $2n+1$ (into odd factors $>1$)": 

$1, 1, 1, 2, 1, 1, 3, 1, 1, 3, 1, 2,\ldots$.

On the other hand, the sequence for odd weights is not yet in \cite{Sl}, namely
 $$\overline{\pp}(2n-1), n\geq 1:\ 2, 5, 5, 13, 5, 19, 5, 34, 13, 19, 5, 65, 5,\ldots $$

Table \ref{tabpp} gives the values of $\overline{\pp}(n,r)$ and $\overline{\pp}(n)$ for certain small $n$ with $0\leq r\leq 5$. Note that $\overline{\pp}(n,r)>0$ if and only if $0\leq r\leq\nu_2(n+1)$.  

The reader may verify some of the values in the last column of the table using the foregoing formulas with the recurrence relation \eqref{off1}.

For example, if $n=23$, then $n+1=2^3 3$ and $s=3$. First, we obtain $f(d)$ for each proper divisor $d$ of 24. One can show directly that $f(1)=f(2)=f(3)=1,\ f(4)=2,\, f(6)=3,\, f(8)=4$ and $f(12)=8$. Then from \eqref{off1}, we have  $f(24)=1+1+1+2+3+4+8=20$. Finally, using \eqref{popform3} we obtain 
$$\overline{\pp}(23) = 4(20)+\left[\binom{3}{3}-(2)^2(4)\right]f(3,1) = 65.$$

\begin{table}[h]
\centering
\footnotesize
\begin{tabular}[t]{|c|c|c|c|c|c|c|c|}
\multicolumn{6}{c}{} \\
\hline
$n/r$&0&1&2&3&4&$\overline{\pp}(n)$\\ \hline
1&1&1& 0& 0& 0& 2\\ \hline
2& 1& 0& 0& 0& 0& 1\\ \hline
3& 2& 2& 1& 0& 0& 5\\ \hline
4& 1& 0& 0& 0& 0& 1\\ \hline
5& 3& 2& 0& 0& 0& 5\\ \hline
6& 1& 0& 0& 0& 0& 1\\ \hline
7& 4& 5& 3& 1& 0& 13\\ \hline
8& 2& 0& 0& 0& 0& 2\\ \hline
9& 3& 2& 0& 0& 0& 5\\ \hline
10& 1& 0& 0& 0& 0& 1\\ \hline
11& 8& 8& 3& 0& 0& 19\\ \hline
12& 1& 0& 0& 0& 0& 1\\ \hline
13& 3& 2& 0& 0& 0& 5\\ \hline
14& 3& 0& 0& 0& 0& 3\\ \hline
15& 8& 12& 9& 4& 1& 34\\ \hline
16& 1& 0& 0& 0& 0& 1\\ \hline
17& 8& 5& 0& 0& 0& 13\\ \hline
18& 1& 0& 0& 0& 0& 1\\ \hline
19& 8& 8& 3& 0& 0& 19\\ \hline
20& 3& 0& 0& 0& 0& 3\\ \hline
21& 3& 2& 0& 0& 0& 5\\ \hline
22& 1& 0& 0& 0& 0& 1\\ \hline
23& 20& 26& 15& 4& 0& 65\\ \hline
24& 2& 0& 0& 0& 0& 2\\ \hline
25& 3& 2& 0& 0& 0& 5\\  \hline
\end{tabular}
$\ $
\begin{tabular}[t]{|c|c|c|c|c|c|c|c|}
\multicolumn{6}{c}{} \\
\hline
$n/r$&0&1&2&3&4&5&$\overline{\pp}(n)$\\ \hline
26& 4& 0& 0& 0& 0& 0& 4\\ \hline
27& 8& 8& 3& 0& 0& 0& 19\\ \hline
28& 1& 0& 0& 0& 0& 0& 1\\ \hline
29& 13& 8& 0& 0& 0& 0& 21\\ \hline
30& 1& 0& 0& 0& 0& 0& 1\\ \hline
31& 16& 28& 25& 14& 5& 1& 89\\ \hline
32& 3& 0& 0& 0& 0& 0& 3\\ \hline
33& 3& 2& 0& 0& 0& 0& 5\\ \hline
34& 3& 0& 0& 0& 0& 0& 3\\ \hline
35& 26& 26& 9& 0& 0& 0& 61\\ \hline
36& 1& 0& 0& 0& 0& 0& 1\\ \hline
37& 3& 2& 0& 0& 0& 0& 5\\ \hline
38& 3& 0& 0& 0& 0& 0& 3\\ \hline
39& 20& 26& 15& 4& 0& 0& 65\\ \hline
40& 1& 0& 0& 0& 0& 0& 1\\ \hline
41& 13& 8& 0& 0& 0& 0& 21\\ \hline
42& 1& 0& 0& 0& 0& 0& 1\\ \hline
43& 8& 8& 3& 0& 0& 0& 19\\ \hline
44& 8& 0& 0& 0& 0& 0& 8\\ \hline
45& 3& 2& 0& 0& 0& 0& 5\\ \hline
46& 1& 0& 0& 0& 0& 0& 1\\ \hline
47& 48& 76& 57& 24& 5& 0& 210\\ \hline
48& 2& 0& 0& 0& 0& 0& 2\\ \hline
49& 8& 5& 0& 0& 0& 0& 13\\ \hline
50& 3& 0& 0& 0& 0& 0& 3\\ \hline
\end{tabular}
\caption{Values of $\overline{\pp}(n,r)$ for $1\leq n\leq 50,\, 0\leq r\leq 5$}\label{tabpp}
\end{table}

\bibliographystyle{amsplain}

\end{document}